\documentclass[12pt,a4paper,twoside]{article}

\newcommand{\pageformat}[6]{\setlength{\hoffset}{-1in}
                  \setlength{\voffset}{-1in}
                  \addtolength{\hoffset}{#5}
                            \addtolength{\voffset}{#6}
                            \setlength{\oddsidemargin}{#1}
                            \setlength{\evensidemargin}{#2}
                            \setlength{\textwidth}{\paperwidth}
                  \addtolength{\textwidth}{-\oddsidemargin}
                  \addtolength{\textwidth}{-\evensidemargin}
                  \addtolength{\textwidth}{-\marginparsep}
                  \addtolength{\textwidth}{-\marginparwidth}
                            \setlength{\topmargin}{#3}
                            \setlength{\textheight}{\paperheight}
                  \addtolength{\textheight}{-\topmargin}
                  \addtolength{\textheight}{-\headheight}
                  \addtolength{\textheight}{-\headsep}
                  \addtolength{\textheight}{-\footskip}
                  \addtolength{\textheight}{-#4}}
\pageformat{2cm}{3cm}{25mm}{25mm}{1pt}{0pt}

\usepackage{ifthen}
\newboolean{article}
    \setboolean{article}{true}
\newboolean{report}
\newboolean{book}
\newboolean{letter}
\newboolean{german}
\newboolean{italian}
\newboolean{nobaselinestretch}
\newboolean{nosectionappendix}
\newboolean{oldtoc}
\newboolean{nosectionequn}
\newboolean{notheorem}

\ifthenelse{\boolean{german}}{
    \usepackage{german}}{}

\usepackage[latin1]{inputenc}

\usepackage{amsmath}
\usepackage{amssymb}
\usepackage[mathscr]{eucal}

\ifthenelse{\boolean{notheorem}}{}{
    \usepackage{theorem}}

\ifthenelse{\boolean{nobaselinestretch}}{}{
    \renewcommand{\baselinestretch}{1.25}}

\newenvironment{env}[2]{\begin{#1}#2\end{#1}}{}
    \newcommand{\beq}[1]{\begin{env}{equation}{#1}}
    \newcommand{\beqn}[1]{\begin{env}{equation*}{#1}}
    \newcommand{\bal}[1]{\begin{env}{align}{#1}}
    \newcommand{\baln}[1]{\begin{env}{align*}{#1}}
    \newcommand{\bga}[1]{\begin{env}{gather}{#1}}
    \newcommand{\bgan}[1]{\begin{env}{gather*}{#1}}
    \newcommand{\bflal}[1]{\begin{env}{flalign}{#1}}
    \newcommand{\bflaln}[1]{\begin{env}{flalign*}{#1}}
    \newcommand{\bmu}[1]{\begin{env}{multline}{#1}}
    \newcommand{\bmun}[1]{\begin{env}{multline*}{#1}}
    \newcommand{\bsp}[1]{\begin{env}{split}{#1}}

    \newcommand{\eeq}{\end{env}}
    \newcommand{\eeqn}{\end{env}}
    \newcommand{\eal}{\end{env}}
    \newcommand{\ealn}{\end{env}}
    \newcommand{\ega}{\end{env}}
    \newcommand{\egan}{\end{env}}
    \newcommand{\eflal}{\end{env}}
    \newcommand{\eflaln}{\end{env}}
    \newcommand{\emu}{\end{env}}
    \newcommand{\emun}{\end{env}}
    \newcommand{\esp}{\end{env}}

\newcommand{\lf}{\vspace{2ex}}

\renewcommand{\bf}[1]{\textbf{#1}}
\renewcommand{\it}[1]{\textit{#1}}

\renewcommand{\sf}[1]{\textsf{#1}}

\renewcommand{\tt}[1]{\texttt{#1}}
\newcommand{\hl}[1]{\bf{\it{#1}}}

\newcommand{\mbf}[1]{\mathbf{#1}}
\newcommand{\msf}[1]{\text{\small$\sf{#1}$}}

\newcommand{\cmc}[1]{\mathcal{#1}}
\newcommand{\eus}[1]{\mathscr{#1}}
\newcommand{\euf}[1]{\mathfrak{#1}}
\newcommand{\bb}[1]{\mathbb{#1}}

\newcommand{\mfootnotesize}[1]{{\setlength{\arraycolsep}{.5ex}\text{\footnotesize$#1$}}}
\newcommand{\mscriptsize}[1]{{\setlength{\arraycolsep}{.3ex}\text{\scriptsize$#1$}}}
\newcommand{\mtiny}[1]{{\setlength{\arraycolsep}{.3ex}\text{\tiny$#1$}}}
\newcommand{\nbd}[1]{$#1$\nobreakdash--}
\newcommand{\ol}[1]{\overline{#1}}

\newcommand{\wt}[1]{\widetilde{#1}}
\newcommand{\wh}[1]{\widehat{#1}}

\newcommand{\om}{\omega}
\newcommand{\Om}{\Omega}

\newcommand{\norm}[1]{\left\lVert#1\right\rVert}

\newcommand{\Babs}[1]{\Bigl\lvert#1\Bigr\rvert}

\newcommand{\snorm}[1]{\norm{\smash{#1}}}

\newcommand{\bfam}[1]{\bigl(#1\bigr)}

\newcommand{\AB}[1]{\langle#1\rangle}

\newcommand{\BAB}[1]{\Bigl\langle#1\Bigr\rangle}
\newcommand{\CB}[1]{\{#1\}}
\newcommand{\bCB}[1]{\bigl\{#1\bigr\}}
\newcommand{\BCB}[1]{\Bigl\{#1\Bigr\}}

\newcommand{\Matrix}[1]{\begin{pmatrix}#1\end{pmatrix}}

\newcommand{\fMatrix}[1]{\mfootnotesize{\Matrix{#1}}}
\newcommand{\sMatrix}[1]{\mscriptsize{\Matrix{#1}}}
\newcommand{\tMatrix}[1]{\mtiny{\Matrix{#1}}}

\newcommand{\set}[2][]{
    \ifthenelse{\equal{#1}{}}{
        \CB{#2}}{
        \CB{#1~|~#2}}}
\newcommand{\bset}[2][]{
    \ifthenelse{\equal{#1}{}}{
        \bCB{#2}}{
        \bCB{#1~|~#2}}}
\newcommand{\Bset}[2][]{
    \ifthenelse{\equal{#1}{}}{
        \BCB{#2}}{
        \BCB{#1~\big|~#2}}}
\newcommand{\zero}{\CB{0}}

\DeclareMathOperator{\ls}{\normalfont\msf{span}}
\DeclareMathOperator{\cls}{\ol{\ls}}

\DeclareMathOperator{\id}{\normalfont\msf{id}}

\newcommand{\C}{\bb{C}}

\newcommand{\N}{\bb{N}}

\newcommand{\R}{\bb{R}}

\newcommand{\cA}{\cmc{A}}
\newcommand{\cB}{\cmc{B}}

\newcommand{\sB}{\eus{B}}

\newcommand{\sK}{\eus{K}}

\newcommand{\ek}{\euf{k}}

\newcommand{\eK}{\euf{K}}
\newcommand{\eL}{\euf{L}}
\newcommand{\eO}{\euf{O}}

\newcommand{\eS}{\euf{S}}
\newcommand{\eT}{\euf{T}}

\newcommand{\U}{\mbf{1}}

\ifthenelse{\boolean{nosectionequn}}{}{
    \numberwithin{equation}{section}
    }

\ifthenelse{\boolean{article}\or\boolean{letter}\or\boolean{nosectionequn}}{
    \setboolean{nosectionappendix}{true}}{}
\ifthenelse{\boolean{nosectionappendix}}{}{
    \renewcommand{\appendix}{
        \chapter*{\appendixname}
        \addcontentsline{toc}{chapter}{\appendixname}
        \renewcommand{\thesection}{\Alph{section}}
        \setcounter{section}{0}}}
   
\ifthenelse{\boolean{report}\or\boolean{book}}{
    }{}

\ifthenelse{\boolean{notheorem}}{}{
        \newcommand{\mnname}{Mathematical note.}
        \newcommand{\enname}{End of the note.}
        \newcommand{\definame}{Definition.}
        \newcommand{\propname}{Proposition.}
        \newcommand{\lemname}{Lemma.}
        \newcommand{\exname}{Example.}
        \newcommand{\exername}{Exercise.}
        \newcommand{\remname}{Remark.}
        \newcommand{\obname}{Observation.}
        \newcommand{\thmname}{Theorem.}
        \newcommand{\corname}{Corollary.}
        \newcommand{\proofname}{Proof.}
    \ifthenelse{\boolean{german}}{
        \renewcommand{\mnname}{Mathematische Notiz.}
        \renewcommand{\enname}{Ende der Notiz.}
        \renewcommand{\exname}{Beispiel.}
        \renewcommand{\exername}{Übung.}
        \renewcommand{\remname}{Bemerkung.}
        \renewcommand{\obname}{Beobachtung.}
        \renewcommand{\thmname}{Satz.}
        \renewcommand{\corname}{Korollar.}
        \renewcommand{\proofname}{Beweis.}}{}
    \ifthenelse{\boolean{italian}}{
        \renewcommand{\mnname}{Nota matematica.}
        \renewcommand{\enname}{Fina della nota.}
        \renewcommand{\definame}{Definizione.}
        \renewcommand{\propname}{Proposizione.}
        \renewcommand{\exname}{Esempio.}
        \renewcommand{\exername}{Esercizio.}
        \renewcommand{\remname}{Nota.}
        \renewcommand{\obname}{Osservazione.}
        \renewcommand{\thmname}{Teorema.}
        \renewcommand{\corname}{Corollario.}
        \renewcommand{\proofname}{Dimostrazione.}

       \renewcommand{\appendixname}{Appendice}

       }{}
    \theoremheaderfont{\normalfont\bfseries}
    \theoremstyle{change}
        \theorembodyfont{\rmfamily}
            \newtheorem{emp}{}[section]
                \newcommand{\bemp}[1][]{
                    \begin{emp}\hskip-\labelsep\bf{#1}\hskip\labelsep}
                \newcommand{\eemp}{\end{emp}}
\newtheorem{itemp}[emp]{}
                \newcommand{\bitemp}[1][]{
                    \begin{itemp}\hskip-\labelsep\bf{#1}\hskip\labelsep\normalfont\itshape}
                \newcommand{\eitemp}{\end{itemp}}
            \newtheorem{mn}[emp]{\mnname}
                \newcommand{\bnm}{\begin{mn}~\begin{quotation}\renewcommand{\baselinestretch}{1}\small\noindent\ignorespaces}
                \newcommand{\enm}{\end{quotation}\hfill\bf{\enname}\end{mn}}
            \newtheorem{ex}[emp]{\exname}
                \newcommand{\bex}{\begin{ex}}
                \newcommand{\eex}{\end{ex}}
            \newtheorem{exer}[emp]{\exername}
                \newcommand{\bexer}{\begin{exer}}
                \newcommand{\eexer}{\end{exer}}
            \newtheorem{defi}[emp]{\definame}
                \newcommand{\bdefi}{\begin{defi}}
                \newcommand{\edefi}{\end{defi}}
            \newtheorem{rem}[emp]{\remname}
                \newcommand{\brem}{\begin{rem}}
                \newcommand{\erem}{\end{rem}}
            \newtheorem{ob}[emp]{\obname}
                \newcommand{\bob}{\begin{ob}}
                \newcommand{\eob}{\end{ob}}
        \theorembodyfont{\normalfont\itshape}
            \newtheorem{thm}[emp]{\thmname}
                \newcommand{\bthm}{\begin{thm}}
                \newcommand{\ethm}{\end{thm}}
            \newtheorem{prop}[emp]{\propname}
                \newcommand{\bprop}{\begin{prop}}
                \newcommand{\eprop}{\end{prop}}
            \newtheorem{cor}[emp]{\corname}
                \newcommand{\bcor}{\begin{cor}}
                \newcommand{\ecor}{\end{cor}}
            \newtheorem{lem}[emp]{\lemname}
                \newcommand{\blem}{\begin{lem}}
                \newcommand{\elem}{\end{lem}}
\newenvironment{empn}[1]{\lf\noindent\bf{#1}\ignorespaces\hskip\labelsep}{\lf}
		\newcommand{\bempn}[1]{\begin{empn}{#1}}
		\newcommand{\eempn}{\end{empn}}
		\newcommand{\bitempn}[1]{\begin{empn}{#1}\normalfont\itshape}
		\newcommand{\eitempn}{\end{empn}}
                \newcommand{\bnmn}{\begin{empn}{\mnname}~\begin{quotation}\renewcommand{\baselinestretch}{1}\small\noindent\ignorespaces}
                \newcommand{\enmn}{\end{quotation}\hfill\bf{\enname}\end{empn}}
		\newcommand{\bexn}{\begin{empn}{\exname}}
		\newcommand{\eexn}{\end{empn}}
		\newcommand{\bexern}{\begin{empn}{\exername}}
		\newcommand{\eexern}{\end{empn}}
		\newcommand{\bdefin}{\begin{empn}{\definame}}
		\newcommand{\edefin}{\end{empn}}
		\newcommand{\bremn}{\begin{empn}{\remname}}
		\newcommand{\eremn}{\end{empn}}
		\newcommand{\bobn}{\begin{empn}{\obname}}
		\newcommand{\eobn}{\end{empn}}

\newcommand{\qedsymbol}{~\rule[-0.35mm]{2mm}{2mm}}
    \newcounter{proof}[emp]
    \newenvironment{Proof}[1]{
        \vspace{1ex}
        \renewcommand{\item}[1][\stepcounter{proof}(\roman{proof})]%
            {##1\hskip\labelsep}
        \noindent\textsc{#1\hskip\labelsep}}{
        \nolinebreak\qedsymbol}
    \newcommand{\proof}[1][\proofname]{
        \begin{Proof}{#1}\ignorespaces}
    \newcommand{\qed}{\end{Proof}}
    \newcommand{\noqed}{
        \renewcommand{\qedsymbol}{}
        \end{Proof}}}
    \ifthenelse{\boolean{italian}}{
        \renewcommand{\proofname}{Dimostrazione.}}{}

\usepackage[varg]{txfonts}

\usepackage[hypertex]{hyperref}

\begin{document}

\bibliographystyle{amsalpha}

\title{The Powers Sum of spatial CPD-semigroups and CP-semigroups}
\author{}
\author{
~\\
Michael Skeide\thanks{This work is supported by research funds of of the Dipartimento S.E.G.e S.\ and Italian MUR under PRIN 2007.}\\\\
{\small\itshape Universit\`a\ degli Studi del Molise}\\
{\small\itshape Dipartimento S.E.G.e S.}\\
{\small\itshape Via de Sanctis}\\
{\small\itshape 86100 Campobasso, Italy}\\
{\small{\itshape E-mail: \tt{skeide@unimol.it}}}\\
{\small{\itshape Homepage: \tt{http://www.math.tu-cottbus.de/INSTITUT/lswas/\_skeide.html}}}\\
}
\date{November 2008}

{
\renewcommand{\baselinestretch}{1}
\maketitle


\vspace{10ex}
\begin{abstract}
\noindent
We define spatial CPD-semigroup and construct their Powers sum. We construct the Powers sum for general spatial CP-semigroups. In both cases, we show that the product system of that Powers sum is the product of the spatial product systems of its factors. We show that on the domain of intersection, pointwise bounded CPD-semigroups on the one side and Schur CP-semigroups on the other, the constructions coincide. This summarizes all known results about Powers sums and generalizes them considerably.
\end{abstract}

}


{\parskip0.5ex plus 0.5ex minus 0.5ex

\section{Introduction}

In the 2002 AMS-Workshop on `Advances in Quantum Dynamics' in Mount Holyoke, Powers described a \it{sum operation} for \it{spatial} \nbd{E_0}semigroups on $\sB(H)$, the algebra of bounded operators on a Hilbert space $H$. The result is a Markov semigroup and Powers asked for the product system of that Markov semigroup in the sense of Bhat \cite{Bha96}, and if that product system coincides or not with the tensor product of the Arveson systems of the \nbd{E_0}semigroups. (By \hl{Arveson system} we shall refer to product systems of Hilbert spaces as introduced by Arveson \cite{Arv89a}, while \hl{product system} refers to the more general situation of Hilbert modules.)

Still during the workshop (see Skeide \cite{Ske03c}) we could show that the Arveson system of the Powers sum is our \it{product} of spatial product systems introduced in \cite{Ske06d} (preprint 2001) immediately for Hilbert modules. In the case of Hilbert spaces, the product is a subsystem of the tensor product. (For modules there is no tensor product of product systems.) Liebscher \cite{Lie00p1} showed that the product may but need not be all of the tensor product. The question, if the subsystem of the tensor product is nevertheless isomorphic to the full tensor product or not, remained open until Powers \cite{Pow04}: It need not.

The Powers sum has been generalized in several directions. Powers \cite{Pow04} generalized it to CP-semigroups that are \it{spatial} in his sense (a sense we consider too narrow). In Bhat, Liebscher and Skeide \cite{BLS02p} we constructed the Powers sum for spatial \nbd{E_0}semigroups on $\sB^a(E)$, the algebra of adjointable operators on a Hilbert module $E$. We also showed that the product system of the sum is our product.

\lf
In Section \ref{sCPDSEC} we introduce \it{spatial} CPD-semigroups and construct their spatial product systems. (This adds several new facts to CPD-semigroups and their GNS-systems as discussed in \cite{BBLS04}. In particular, like discussed in \cite{BLS08p} for spatial CP-semigroups, the spatial product system of a CPD-semigroup, in the \nbd{C^*}case, may be bigger than the GNS-system.) In Section \ref{CPDPsSEC} we construct a Powers sum for them, which is a spatial CPD-semigroup, too. We show that the product system of the sum is our product of the spatial product systems of the constituents. In Section \ref{CPPsSEC} we introduce a Powers sum for spatial strict CP-semigroups acting on (not necessarily equal) $\sB^a(E)$s, and show also here that their the spatial product systems of the sum is the product of the spatial product systems of the constituents. Both sorts of Powers sums include Powers construction \cite{Pow04} (adding to \cite{Pow04} the identification of the Arveson system of the sum) and generalize it considerably. Our second Powers sum for CP-semigroups on $\sB^a(E)$ includes and generalizes \cite{BLS02p} and furnishes the case treated there, \nbd{E_0}semigroup, with a more transparent proof. Finally, in Section \ref{vsSEC} we show that the subclass of pointwise bounded CPD-semigroups and the subclass of Schur CP-semigroups are two sides of the same coin.

The discussion is mainly for \nbd{C^*}algebras and modules. For the reasons explained in Section \ref{sCPDSEC}, this case is more peculiar. With few modifications, also explained in Section \ref{sCPDSEC}, the case of von Neumann algebras and modules is always included, usually, with simplified proofs.

It would be interesting to follow the story in the historical order. But for this we would have to introduce a lot of terminology, needed just to describe the known results, before we came to new ones. We prefer, therefore, to start immediately with the discussion of spatial CPD-semigroups and their product systems, followed the definition of their Powers sum and the identification of its product system. Only then we explain how this specializes to Powers results.
}

\section{Spatial CPD-semigroups and their product systems}\label{sCPDSEC}

Let $S$ denote a set. We shall consider kernels $\eK\colon S\times S\rightarrow\sB(\cA,\cB),(s,s')\mapsto\eK^{s,s'}$ with values in the bounded maps between two \nbd{C^*}algebras $\cA$ and $\cB$. In the case of von Neumann algebras we shall require the maps $\eK^{s,s'}$ to be \nbd{\sigma}weak. Following the definition in Barreto, Bhat, Liebscher and Skeide \cite{BBLS04}, we say a kernel is \hl{completely positive definite} (\hl{CPD}), if
\beq{\label{CPDdef}
\sum_{i,j}b_i^*\eK^{s_i,s_j}(a_i^*a_j)b_j
~\ge~
0
}\eeq
for all choices of finitely many elements $s_i\in S,a_i\in\cA,b_i\in\cB$. A typical example of a CPD-kernel is given by
\beqn{
\eK^{s,s'}
~:=~
\AB{\xi^s,\bullet\xi^{s'}}
}\eeqn
for a family $\bfam{\xi^s}_{s\in S}$ of elements in a correspondence $E$ from $\cA$ to $\cB$. In fact, if $\cA$ and $\cB$ are unital, then every CPD-kernel can be recovered in that way. by its \hl{Kolmogorov decomposition}. If we require that $E$ is generated by $\bfam{\xi^s}_{s\in S}$ as a correspondence, then the pair $(E,\bfam{\xi^s}_{s\in S})$ is unique up to bilinear unitary equivalence. We refer to it as the \hl{Kolmogorov decomposition} and to $E$ as the \hl{Kolmogorov correspondence} of $\eK$.

\brem\label{kerrem}
If $S$ has exactly one element, then the CPD-kernels on $S$ are precisely the CP-maps, and the Kolmogorov decomposition is Paschke's GNS-construction for CP-maps; see \cite{Pas73}. If $S=\CB{1,\ldots,n}$, then the CPD-kernels on $S$ can be identified with \hl{Schur} CP-maps from $M_n(\cA)$ to $M_n(\cB)$ that act matrix elementwise as $a_{ij}\mapsto\eK^{i,j}(a_{i,j})$. If $\cA=\C$, then by $(s,s')\mapsto\eK^{s,s'}(1)$ we establish a one-to-one correspondence with positive definite \nbd{\cB}valued kernels. If also $\cB=\C$, then we get back the usual Kolmogorov decomposition for \nbd{\C}valued kernels.
\erem

\brem
Even if $\cA$ and $\cB$ are nonunital, we get a correspondence $E$ from $\cA$ to $\cB$ and elements $\xi^{a,s}$ such that $\eK^{s,s'}(a^*a')=\AB{\xi^{a,s},\xi^{a',s'}}$. But it is, in general, impossible to obtain suitable elements $\xi^s$. It is possible to unitalize the kernel to the unitalizations $\wt{\cA}$ and $\wt{\cB}$ by the unitalization procedure in Skeide \cite{Ske08} or, if all $\eK^{s,s'}$ are strict, to the multiplier algebras.
\erem

A \hl{CPD-semigroup} is a family $\eT=\bfam{\eT_t}_{t\in\R_+}$ of \nbd{\sB(\cB)}valued CPD-kernels on $S$ such that for each $s,s'\in S$ the maps $\eT^{s,s'}_t$ form a semigroup on $\cB$. If all these semigroups are continuous in a certain topology, then we say the CPD-semigroup is continuous in that topology.

Like CPD-kernels are related to correspondences via Kolmogorov decomposition, CPD-semigroups are related to product systems of correspondences via a GNS-type construction. Following Bhat and Skeide \cite{BhSk00}, a \hl{product system} is a family $E^\odot=\bfam{E_t}_{t\in\R_+}$ of correspondences over $\cB$ such that $E_0=\cB$ with a family of bilinear unitaries $u_{s,t}\colon E_s\odot E_t\rightarrow E_{s+t}$ such that the product $x_sy_t:=u_{s,t}(x_s\odot y_t)$ is associative, and such that $u_{0,t}$ and $u_{t,0}$ are the canonical identifications. If $\cB$ is unital, a \hl{unit} for $E^\odot$ is a family $\xi^\odot=\bfam{\xi_t}_{t\in\R_+}$ of elements $\xi_t\in E_t$ such that $\xi_0=\U$ and such that $\xi_s\xi_t=\xi_{s+t}$. If $\bfam{{\xi^s}^\odot}_{s\in S}$ is a family of units for $E^\odot$, then by the definition of the internal tensor product $E_s\odot E_t$ it follows that the maps
\beqn{
\eT^{s,s'}_t
~:=~
\AB{\xi^s_t,\bullet\xi^{s'}_t}
}\eeqn
form a semigroup. Clearly, the family of kernels $(s,s')\mapsto\eT^{s,s'}_t$ forms a CPD-semigroup. By \cite{BBLS04}, every CPD-semigroup for unital $\cB$ arises in that way. If $E^\odot$ is generated as a product system by the family of units $\bfam{{\xi^s}^\odot}_{s\in S}$, then the pair $(E^\odot,\bfam{{\xi^s}^\odot}_{s\in S})$ is determined up to suitable isomorphism. We refer to it as the \hl{GNS-construction} and to $E^\odot$ as the \hl{GNS-system} of $\eT$.

\brem
If $S$ has one element, then the CPD-semigroups on $S$ are precisely the CP-semigroups and the GNS-construction is that from \cite{BhSk00}.
\erem

\bob\label{gensetob}
If $E^\odot$ is the GNS-system of a CPD-semigroup $\eT$ and $\bfam{{\xi^s}^\odot}_{s\in S}$ the generating family of units, then
\beqn{
E_t
~=~
\cls\Bset[\,b_n\xi^{s_n}_{t_n}\odot\ldots\odot b_1\xi^{s_1}_{t_1}b_0]{n\in\N;s_1,\ldots,s_n\in S;b_0,\ldots,b_n\in\cB;t_1+\ldots+t_n=t}.
}\eeqn
This is very important to identify, later on, the spatial product system of the Powers sum of spatial CPD-semigroups on $\cB$ or spatial CP-semigroups on $\sB^a(E)$.
\eob

The definitions and results repeated so far, were for (unital) \nbd{C^*}algebras and \nbd{C^*}modules or correspondences. They modify easily to von Neumann algebras, modules, and correspondences, if we: 1.) Require maps on or between von Neumann algebras \nbd{\sigma}weak. 2.) Replace the tensor product by its strongly closed version. Without further mention, we assume these \bf{conventions} when we speak about the von Neumann case.

The following definitions  generalizes Arveson's \cite{Arv97a} for (normal) CP-semigroups on $\sB(G)$ ($G$ some Hilbert space). It is new, except for the domination of CPD-semigroups from \cite{BBLS04}. The version for CP-semigroups on \nbd{C^*}algebras is from Bhat, Liebscher and Skeide \cite{BLS08p}; that for von Neumann algebras $\cB\subset\sB(G)$ from Skeide \cite{Ske08p1}. 

\bdefi
Let $\cB$  be unital \nbd{C^*}algebra (a von Neumann algebra) and let $S$ be a set.

A CPD-semigroup $\eT$ on $S$ with values in $\sB(\cB)$ \hl{dominates} another $\eS$, if the kernels $\eT_t-\eS_t$ are CPD for all $t\in\R_+$. In this situation we write $\eT\ge\eS$.

A CPD-semigroup $\eS$ is \hl{elementary}, if it has the form $\eS^{s,s'}_t={c^s_t}^*\bullet c^{s'}_t$ for a family $\bfam{c^s}_{s\in S}$ of (strongly) continuous semigroups $c^s=\bfam{c^s_t}_{t\in\R_+}$ in $\cB$.

A \hl{unit} for a CPD-semigroup $\eT$ is an elementary CPD-semigroup $\eS$ such that $\eT\ge\eS$.

A CPD-semigroup $\eT$ is \hl{spatial}, if it admits a unit. If we wish to emphasize the choice of the unit, we will also speak of the pair $(\eT,\eS)$ as spatial CPD-semigroup.
\edefi

\brem
As pointed out in \cite{BLS08p}, semigroups of elements $c_t$ in a \bf{unital} \nbd{C^*}algebra that are continuous in any of the natural topologies, are uniformly continuous automatically. Indeed, if $c_t$ is weakly continuous, then the semigroup $b\mapsto bc_t$ of maps in $\sB(\cB)$ it is weakly continuous and, therefore, strongly continuous. In particular, the family $c_t=\U c_t\in\cB$ is norm continuous. (If $\cB$ is nonunital, then it has no sense to speak of a semigroup in $\cB$ indexed by $t\ge0$, but only $t>0$.)  The strong topology of a von Neumann algebra $\cB\subset\sB(G)$ is much weaker and allows for semigroups with unbounded generator.

Note, too, that spatiality without continuity conditions on the unit $\eS$ is a trivial issue. In fact, the \hl{zero-semigroup} $\eO$ defined by $\eO^{s,s'}_t=0$ for all $s,s\in S$ and all $t>0$ would be a unit for every CPD-semigroup on $S$.
\erem

In the sequel, \hl{strongly continuous} for a semigroup $T$ of maps on a \nbd{C^*}algebra $\cB$ that is not represented as an algebra of operators on a Hilbert space or, more generally, on a Hilbert module, means that $t\mapsto T_t(b)$ is norm continuous for every $b\in\cB$. If the $T$ acts on a von Neumann algebra $\cB\subset\sB(G)$, then we mean that $t\mapsto T_t(b)g$ is norm continuous for all $b\in\cB,g\in G$. The same convention applies to semigroups acting on $\sB^a(E)$. We see in a minute that \bf{spatial} strongly continuous CPD-semigroup on an abstract \nbd{C^*}algebra are even uniformly continuous. Since in these notes we are interested only in spatial CPD-semigroups, we will, usually, use strongly continuous only when we speak about operator algebras $\cB\subset\sB(G)$ or $\sB^a(E)$.

\bemp[Definition \cite{Ske06d}.~]
A \hl{spatial} product system is a pair $(E^\odot,\om^\odot)$ consisting of a product system and a central unital reference unit $\om^\odot$ (that is, all $\om_t$ commute with all $b\in\cB$ and all $\om_t$ are unit vectors in the sense that $\AB{\om_t,\om_t}=\U$).
\eemp

In general, we will call a product system spatial, if it has central unital units. Note, however, that the spatial structure may depend on the choice of the reference unit. Spatial product systems and their product (see Section \ref{CPDPsSEC}) have been introduced in Skeide \cite{Ske06d}. They form a subcategory of product systems which behaves best in analogy with Arveson's classification scheme for Arveson systems. (There is an \it{index} for spatial product systems and the index behaves additively under the \it{product} of spatial product systems; see \cite{Ske06d}. A tensor product of product systems does, in general, not exist.)

\bthm\label{spatCPDthm}
For a strongly continuous CPD-semigroup $\eT$ on $S$ with values in $\sB(\cB)$ for a unital \nbd{C^*}algebra $\cB$ the following conditions are equivalent:
\begin{enumerate}
\item\label{sCPD1}
$\eT$ is spatial.

\item\label{sCPD2}
The (continuous) GNS-system of $\eT$ embeds into a (continuous) spatial product system. In particular, $\eT$ is uniformly continuous.

\item\label{sCPD3}
$\eT$ has a \hl{Christensen-Evans generator}, that is, $\eL^{s,s'}:=\frac{d\eT^{s,s'}}{dt}\big|_{t=0}$ exists for all $s,s'\in S$ and there are a CPD-kernel $\eL_0$ and elements $\beta_s\in\cB$ such that
\beqn{
\eL^{s,s'}(b)
~=~
\eL_0^{s,s'}(b)+b\beta_{s'}+\beta_s^*b.
}\eeqn

\end{enumerate}
\ethm

\proof
The proof is very much like the proofs of \cite[Theorem 3.4 and Corollary 3.7]{BLS08p}, just more indices. (The number of indices is $\#S+1$.) In so far, we explain only the construction of the extended CPD-semigroup on $S_0=S\cup\zero$, and we say a word on what continuous product systems means.

We start with the latter. If a CPD-semigroup is strongly continuous, then its product system is continuous in the sense of Skeide \cite{Ske03b}. By \cite[Theorem 7.7]{Ske03b}, if a continuous product system has a single unit $\xi^\odot$ such that the CP-semigroup $\AB{\xi_t,\bullet\xi_t}$ is uniformly continuous, then all semigroups $\AB{\xi_t,\bullet\xi'_t}$ are uniformly continuous. And the reference unit $\om^\odot$ generates the trivial semigroup which is uniformly continuous.

The basic observation for constructing the spatial product system into which the GNS-system embeds, is the following. Let $\eS$ be an elementary CPD-semigroup on $S$ generated by semigroups $c^s$ in $\cB$, and suppose that $\eT$ dominates $\eS$. Then the semigroup $\wh{\eT}$ on $S_0$ defined by setting
\baln{
\wh{\eT}_t^{s,s'}
&
~:=~
\eT_t^{s,s'},
&
\wh{\eT}_t^{0,s}
&
~:=~
\bullet c_t,
&
\wh{\eT}_t^{s,0}
&
~:=~
c_t^*\bullet,
}\ealn
is CPD. (It can be written as the sum of the extension of $\eT_t-\eS_t$ from $S$ to $S_0$ by $0$, and a suitable elementary CPD-semigroup on $S_0$; see \cite{BLS08p}.) Clearly, the GNS-system of $\wh{\eT}$ is spatial (the unit ${\xi^0}^\odot$ is central and unital), and it contains the GNS-system of $\eT$, see \cite{BLS08p} for details.\qed

\brem
In general, the generators of uniformly continuous CPD-semi\-groups with values in $\sB(\cB)$ are precisely the \hl{conditionally completely positive definite} (\hl{CCPD}) \nbd{\sB(\cB)}valued kernels (that is, the kernel fulfills \eqref{CPDdef} under the \it{condition} that $\sum_ia_ib_i=0$); see \cite{BBLS04}. Like for CP-semigroups on a \nbd{C^*}algebra, boundedness of the generator is not sufficient for that is has Christensen-Evans form.
\erem

It is an open problem, whether or not the spatial extension of the GNS-system of $\eT$ constructed in the proof of Theorem \ref{spatCPDthm} depends on the choice of the unit $\eS$. On the other hand, it is easy to see that it does not depend on the choice of the implementing semigroup $c$. (In fact, two semigroups in $\cB$ implementing the same elementary CP-semigroup on $\cB$, can differ at most by a unitary semigroup in the center of $\cB$. Using this, it is more or less obvious to see that GNS-system of the extended kernel $\wh{\eT}$ does not change und such a variation.) Henceforth, we call it \hl{the spatial extension} of the GNS-system \hl{based} on the unit $\eS$. Fortunately, the Powers sum of CPD-semigroups depends from the beginning on the choice of units. So, it is not tragic if also the spatial extension of their GNS-systems should depend on that choice.

For von Neumann algebras the situation is much better:

\bthm
For a strongly continuous CPD-semigroup $\eT$ on $S$ with values in $\sB(\cB)$ for a von Neumann algebra $\cB\subset\sB(G)$ the following conditions are equivalent:
\begin{enumerate}
\item
$\eT$ is spatial.

\item
The GNS-system of $\eT$ is spatial.

\end{enumerate}
\ethm

\proof
The proof is the same as in Skeide \cite{Ske08p1} for CP-semigroups. It cannot be reproduced here for reasons of space. (Very roughly, the idea is that for von Neumann algebras there is an order isomorphism from the partially ordered set of positive contraction endomorphisms of the GNS-system of a CPD-semigroup $\eT$ to the partially ordered set of CPD-semigroups dominated by $\eT$; see \cite{BBLS04}. And the range of the positive contraction morphism corresponding to a unit is just the one-dimensional product system $\bfam{\cB}_{t\in\R_+}$, which contains the central unital unit $\bfam{\U}_{t\in\R_+}$.)\qed

\section{The Powers sum of spatial CPD-semigroups}\label{CPDPsSEC}

Let $\eT^1$ and $\eT^2$ be spatial CPD-semigroups on sets $S^1$ and $S^2$, respectively, with values in $\sB(\cB)$. Choose units $\eS^1$ and $\eS^2$ for them implemented by semigroups $c^1$ and $c^2$, respectively, in $\cB$. Define a semigroup $\eT:=\eT^1\boxplus\eT^2$ on $S:=S^1\uplus S^2$ (disjoint union), by setting
\baln{
(\eT^1\boxplus\eT^2)^{s,s'}
&
~:=~
(\eT^i)^{s,s'}
&&
~~
(i=1,2;s,s'\in S^i),
\\
(\eT^1\boxplus\eT^2)^{s_i,s_j}
&
~:=~
{(c^i)^{s_i}}^*\bullet (c^j)^{s_j}
&&
(i\ne j,s_i\in S^i,s_j\in S^j).
}\ealn
Observe that each $\eS^i$ is itself a spatial CPD-semigroup with unit $\eS^i$. Therefore, the definition applies also to $\eS:=\eS^1\boxplus\eS^2$.

\bthm\label{PsCPDthm}
$\eT$ is a spatial CPD-semigroup with unit $\eS$. Clearly, $\eT$ is (strongly) continuous if and only if each $\eT^i$ is (strongly) continuous.
\ethm

\proof
We shall show $\eT\ge\eS$. This settles both that $\eT_t$ is CPD (as sum of the CPD-kernels $\eT_t-\eS_t$ and $\eS_t$) and that $\eS$ is a unit for $\eT$. We find
\baln{
(\eT_t-\eS_t)^{s,s'}
&
~=~
(\eT^i_t-\eS^i_t)^{s,s'}
&&
~~
(i=1,2;s,s'\in S^i),
\\
(\eT_t-\eS_t)^{s_i,s_j}
&
~=~
0
&&
(i\ne j,s_i\in S^i,s_j\in S^j).
}\ealn
Since each $\eT^i_t-\eS^i_t$ is CPD on $S^i$ and since all mixing terms $s^i\in S^i,s^j\in S^j$ $(i\ne j)$ disappear, this shows that $\eT_t-\eS_t$ is CPD.\qed

\bdefi
We refer to $(\eT,\eS):=(\eT^1,\eS^1)\boxplus(\eT^2,\eS^2)$ as the \hl{Powers sum} of $(\eT^1,\eS^1)$ and $(\eT^2,\eS^2)$.
\edefi

We now wish to identify the spatial extension of the GNS-system of $(\eT,\eS)$ as the product of the spatial extensions of the GNS-systems of $(\eT^i,\eS^i)$. To that goal we repeat the characterization in Skeide \cite{Ske06d} of the product in terms of a universal property.

\bitemp[Theorem and Definition {\cite[Theorem 5.1 and Definition 5.2]{Ske06d}}.~]\label{PRODthm}
Let $({E^1}^\odot,{\om^1}^\odot)$ and $({E^2}^\odot,{\om^2}^\odot)$ denote spatial product systems. Then there exists a spatial product system $(F^\odot,\om^\odot)$ fulfilling the following properties:

\begin{enumerate}
\item\label{psPS1}
$(F^\odot,\om^\odot)$ contains $({E^i}^\odot,{\om^i}^\odot)$ $(i=1,2)$ as spatial subsystems and is generated by them, that is, $F_t$ is spanned by expressions like
\beqn{
x_{t_n}^n\odot\ldots\odot x_{t_1}^1,
}\eeqn
$n\in\N$, $t_1+\ldots+t_n=t$, and $x_{t_j}^j\in E^1_{t_j}\cup E^2_{t_j}$.

\item\label{psPS2}
The inner product of members $x^1\in E^1_t\subset F$ and $x^2\in E^2_t\subset F$ is given by
\beqn{
\AB{x^1,x^2}
~=~
\AB{x^1,\om^1_t}\AB{\om^2_t,x^2}.
}\eeqn
\end{enumerate}
Moreover, every spatial product system fulfilling these properties is canonically isomorphic to $(F^\odot,\om^\odot)$.

We call $(F^\odot,\om^\odot)$ the \hl{product} of the factors $({E^i}^\odot,{\om^i}^\odot)$ and we denote it by $(\bfam{{E^1}\circledcirc{E^2}}^\odot,\om^\odot)$.
\eitemp

Note that, by Property \ref{psPS2}, in the product the two reference units ${\om^i}^\odot$ of the factors get identified with the reference unit $\om^\odot$.

\bthm
The spatial extension of the GNS-system of $(\eT,\eS):=(\eT^1,\eS^1)\boxplus(\eT^2,\eS^2)$ is (spatially) isomorphic to the product of the spatial extensions of the GNS-systems of $(\eT^1,\eS^1)$ and $(\eT^2,\eS^2)$.
\ethm

\proof
Recall that \it{spatially} isomorphic means that the isomorphism identifies also the reference units.

In Observation \ref{gensetob} we indicated a spanning subset of the GNS-system of a CPD-semigroup. We apply this to the GNS-system of the spatial extension of $(\eT,\eS)$. Observe that the pieces of units $\xi^j_{t_j}$ come either from the GNS-system of $\eT^1$ or from the GNS-system of $\eT^2$ or from the component $0$ in $S\cup\zero$, that is from the reference unit of the GNS-system of the spatial extension of $(\eT,\eS)$. One easily verifies that the inner product does not change, if instead we replace that reference unit with one (no matter which) of the reference units of the spatial extensions of the GNS-systems of one of the factors. This shows that the spatial extension of the GNS-system of $(\eT,\eS)$ contains the spatial extensions of the GNS-systems of the factors as subsets and is generated by them, as required in Property \ref{psPS1} of Theorem \ref{PRODthm}. It is also easy to check that the inner products of elements from different factors are those required by Property \ref{psPS2} of Theorem \ref{PRODthm}.\qed

\brem
Notation and formulation of the results is for the \nbd{C^*}case. But this case is the more complicated because the GNS-system of spatial CPD-semigroup need not be spatial. With the standard topological conventions we applied in the preceding section, all statements (some of them in a simpler form) remain valid in the von Neumann case.
\erem

\brem\label{Afamrem}
It is clear that both product of spatial product system and Powers sum of spatial CPD-semigroups may be carried out for families of products systems or CPD-semigroup indexed by arbitrary sets. For spatial product systems this is mentioned in \cite[Remark 5.7]{Ske06d}. For spatial CPD-semigroups this means that for a family $\bfam{\eT^\alpha,\eS^\alpha}_{\alpha\in A}$ there is a Powers sum $\text{\Large$\boxplus$}_{\alpha\in A}\eT^\alpha$. Of course, also the statement of the theorem remains valid for such families: The spatial extension of the GNS-system of the Powers sum is the product of the spatial extensions of the GNS-systems of the constituents.
\erem

\section{Some special cases}\label{exSEC}

In this section we discuss some examples. We have a look how Schur semigroups of positive definite kernels are included. In Remark \ref{nonPDrem}, we explain why such semigroups do not make sense in a noncommutative context, underlining CPD-semigroups as the correct generalization. In Example \ref{SchurCPex} we discuss how the case of finite sets can be described equivalently in terms of \it{Schur CP-semigroup} on matrix algebras. This includes all case discussed by Powers for the case $\cB=\sB(G)$ and even generalizes them in that case. It does \bf{not} cover the case studied in Bhat, Liebscher and Skeide \cite{BLS02p}, but gives a hint what to do in the following section.

\bemp[Example. Semigroups of positive definite kernels.]
If $\cB=\C$, we recover the notion of \hl{positive definite} (\hl{PD}) \nbd{\C}valued kernels and their \hl{Schur semigroups}. In fact, a map $\eK^{s,s'}$ on $\C$ is determined by the value $\ek^{s,s'}:=\eK^{s,s'}(1)$ and a \nbd{\C}valued kernel $\ek$ on $S$ defines a \nbd{\sB(\C)}valued kernel $\eK$ by setting $\eK^{s,s'}(z):=\ek^{s,s'}z$. Clearly, $\eK$ is CPD if and only if $\ek$ is PD.

The \hl{Schur product} of two \nbd{\C}valued kernels on $S$ is simply the pointwise product of functions on $S\times S$. Clearly, the Schur product is reflected by the composition of the corresponding \nbd{\sB(\C)}valued kernels. Everything we know about CPD-semigroups has, thus, an immediate interpretation in terms of PD-semigroups: 1.) The Schur product preserves PD. 2.) PD-semigroups have a product systems of Hilbert spaces (that is, an Arveson system) as GNS-system. 3.) This Arveson system is generated by its units an, therefore Fock. This fact has already been noted by Parthasarathy and Schmidt \cite{PaSchm72}. They applied this knowledge to the PD-semigroup $\ek^{s,s'}_t:=\int e^{i(s-s')}\mu_t(dt)$ on $\R$ that arises from the convolution semigroup $\bfam{\mu_t}_{t\in\R_+}$ of distributions of a Lévy process, which enabled them to represent every Lévy process on the Fock space.

The product of Arveson systems of Fock type, so-called \hl{type I} Arveson systems, is simply their tensor product. (This need not be so for non-type I spatial Arveson systems, so-called \hl{type II} Arveson systems.) Tensor products of units in the factors give rise to units in the tensor product, and every unit in the tensor product arises in that way. In order to understand the PD-semigroup on $S^1\uplus S^2$ it is better to assume that in each factor a reference unit has been distinguished, that corresponds to $0\in S^i\cup\zero$. The set $S^1\uplus S^2$ is, identified with the subset $(S^1\times\zero)\cup(\zero\times S^2)$ of the generating set $(S^1\cup\zero)\times(\zero\cup S^2)$ of units in the product. If both kernels come from Lévy processes, then the product simply describes the two processes as a pair of independent Lévy processes (or a two-dimensional Lévy process) on the same probability space (the product space).

The structure of the units as products of the units of the factors remains valid for the product of arbitrary spatial product systems; see \cite[Theorem 5.6]{Ske06d}. That the disjoint union of two sets $S^1$ and $S^2$ is, by the very definition of \it{disjoint union}, a subset of product of \it{dotted} dotted sets $S^1\cup\zero$ and $\zero\cup S^2$, is reflected in a striking way by the structure of the set of units in the product.
\eemp

\brem\label{nonPDrem}
It is natural to ask for PD-semigroups of \nbd{\cB}valued kernels, sitting somehow in between \nbd{\C}valued PD-kernels and general CPD-kernels. However, among PD-kernels there is no whatsoever product operation that would respect the PD-condition, as soon as the algebra $\cB$ is noncommutative.

This is a central thread of positivity in a noncommutative setting: If we wish to compose positive noncommuting things in a positivity preserving way then \it{composition} must be \it{composition of maps} on the \nbd{*}algebra. Almost never it can be based on multiplication of positive elements in a \nbd{*}algebra. 

The basic feature of positive things is that they possess whatsoever kind of square root. A positive element $a$ in a \nbd{C^*}algebra can be written as $b^*b$. But, if there is another one $a'=b'^*b'$, then $aa'=b^*bb'^*b'$ is almost never positive. An element that is positive, is $b^*b'^*b'b$. However, this element depends manifestly on the choice of the square roots $b$ and $b'$.

A way out is to consider, from the beginning, the map $T=b^*\bullet b$ instead of $a=T(\U)$. The knowledge of $T(\U)$ is only rarely a suitable substitute for the whole map $T$. But, once we have that map, we my compose it with $T'$, and, in fact, we get $T\circ T'=(b'b)^*\bullet(b'b)$.

Going one step further to CP-type maps (for instance, CPD-kernels), one sees that the related GNS-constructions play the role of the square roots which may be \it{multiplied}. The \it{multiplication} is simply the tensor product of the associated GNS-correspondences; see \cite[Observation 2.17]{BhSk00} or \cite[Observation 3.4.3]{BBLS04}. To say it more provocantly: GNS-systems are square roots of the CPD-semigroups they stem from.
\erem

\bemp[Example. Schur CP-semigroups on $M_n(\cB)$.~]\label{SchurCPex}
Recall that the case of a CPD-semigroup on a one-point set $S$, is precisely the case of a CP-semigroup. More generally, a CPD-semigroup on an \nbd{n}point set $S=\CB{1,\ldots,n}$ gives rise to a CP-semigroup $T^n$ on $M_n(\cB)$, by setting
\beqn{
\bfam{T^n_t(A)}_{ij}
~=~
\eT^{i,j}_t(a_{ij}).
}\eeqn
Clearly, we do not obtain all CP-semigroups on $M_n(\cB)$ in that way. In fact, $T^n$ is a \hl{Schur semigroup}, that is, it acts matrix element wise on the matrix $A=\bfam{a_{ij}}$. So, what we really have, is a one-to-one correspondence between CPD-semigroups on a fixed \nbd{n}point set $S$ and Schur CP-semigroups on $M_n(\cB)$.

The elementary CP-semigroups on $M_n(\cB)$ which are also Schur semigroups are precisely those that are generated are generated $S^n_t=C_t^*\bullet C_t$ by semigroups $C=\bfam{C_t}_{t\in\R_+}$ in $M_n(\cB)$ with diagonal matrices $C_t\in M_n(\cB)$. It is easy to check that Schur CP-semigroup is spatial if and only if the corresponding CPD-semigroup $\eT$ on the \nbd{n}point set $S$ is spatial. (The entries of the diagonal generate the elementary CPD-semigroup $\eS$ dominated by $\eT$, and \it{vice versa}.)

We see that there is a Powers sum of spatial Schur CP-semigroups $T^{n_i}$ on $M_{n_i}(\cB)$ that provides a Schur CP-semigroup $T^{n_1}\boxplus T^{n_2}$ acting on $M_{n_1+n_2}(\cB)$.

The special case $n_1=n_2=1$ (semigroups on $\cB$ having a semigroup on $M_2(\cB)$ as sum) includes all cases discussed by Powers with $\cB=\sB(G)$. (In that case, when $G$ is infinite-dimensional and separable, $M_2(\sB(G))=\sB(G\oplus G)$ and $G\oplus G\cong G$. This abuse, mixing a true binary operation among semigroups on $\sB(G)$ with a binary operation among conjugacy classes, is quite common. For more general algebras where, usually, $M_2(\cB)\ncong\cB$, this is no longer possible.) In the Mount Holyoke meeting, Powers proposed the case when the CP-semigroups are spatial semigroups of unital endomorphisms where the units are isometric. In \cite{Pow04} he generalized to spatial CP-semigroups, but still with isometric units. He even calls these CP-semigroups \it{spatial}. But, we think that Arveson's (much) wider definition in \cite{Arv97a} is the adequate one, and our discussion does already extend the Powers sum to that case.
\eemp

In the preceding example we did not say a word about the product systems of the involved semigroups. In fact, the product system of $T^n$ consists correspondences over $M_n(\cB)$, while the product system of the corresponding CPD-semigroup $\eT$ consists of correspondences over $\cB$. For spatial product systems of correspondences over the \it{same} algebra $\cB$, there is the product of spatial product systems. But, the algebras $M_n(\cB)$ may be nonisomorphic for different $n$. How are the product system of $T^n$ and of $\eT$ related, so that the product operation of the product system in the CPD-picture can be applied? Also the question, whether Example \ref{SchurCPex} can be generalized to arbitrary index sets, is interesting. We answer these and other questions in the more general setting of the following section.

\section{The Powers sum of CP-semigroups on $\sB ^a(E)$}\label{CPPsSEC}

Observe that $M_n(\cB)=\sB^a(\cB^n)$. In this section we will replace $\cB^n$ with a general \hl{full} Hilbert \nbd{\cB}module (that is, $\cls\AB{E,E}=\cB$, respectively, $\cls^s\AB{E,E}=\sB$ in the von Neumann case). But, if we do so, then the terminology \it{Schur CP-semigroup} has no longer sense. (This is something which has sense only with respect to an ONB or, possibly, a quasi ONB.) On the other hand, for the \nbd{C^*}case in this setting it is indispensable that we require the CP-semigroups $T$ on $\sB^a(E)$ to be \hl{strict}, that is, each $T_t$ is \nbd{*}strongly continuous on bounded subsets. (In the von Neumann case our standard hypothesis, normality, is sufficient.) The result that on $M_n(\cB)$ the strict topology coincides with the norm topology ($\cB$ is assumed unital!), is standard. This is why, in Example \ref{SchurCPex}, we did not worry about strictness.

Before we study spatial CP-semigroups on $\sB^a(E)$, we first repeat some results  from Bhat, Liebscher, and Skeide \cite{BLS02p} about general strict CP-semigroups on $\sB^a(E)$. We will also derive some new results on spatiality of such semigroups.

In \cite{BLS02p} we showed that the product system $F^\odot$ of a strict CP-semigroup $T$ on $\sB^a(E)$ (consisting of correspondences over $\sB^a(E)$!) may be \it{transformed} into a product system $E^\odot$ consisting of correspondences over $\cB$ in the following way: For each $F_t$ define the \nbd{\cB}cor\-re\-spond\-ence $E_t:=E^*\odot F_t\odot E$, where $E^*$ is a correspondence from $\cB$ to $\sB^a(E)$ with inner product $\AB{x^*,y^*}:=xy^*$ (the \hl{rank-one operator} that maps $z$ to $x\AB{y,z}$) and bimodule action $bx^*a:=(a^*xb^*)^*$. (Note that both tensor products are over $\sB^a(E)$. Note, too, that $E^*\odot E=\cB$ via $x^*\odot y=\AB{x,y}$ and $E\odot E^*=\sK(E)$, the \nbd{C^*}algebra of \hl{compact} operators on $E$, via $x\odot y^*=xy^*$. Since all $T_t$ are strict, the left action of $\sK(E)$ on $F_t$ is nondegenerate.) The $E_t$ form a product system $E^\odot$ via
\beqn{
E_s\odot E_t
~=~
E^*\odot F_s\odot E\odot E^*\odot F_t\odot E
~=~
E^*\odot F_s\odot F_t\odot E
~=~
E^*\odot F_{s+t}\odot E
~=~
E_{s+t}.
}\eeqn
Note that in making disappear the part $E\odot E^*=\sK(E)$ in the middle, we did use strictness of the left action on $F_t$. The right action of $\sK(E)$ on $F_s$ will rarely be nondegenerate:

\brem
The transition $F^\odot\rightsquigarrow E^\odot$ is very close to an operation of \it{Morita equivalence} of product systems, as defined in Skeide \cite{Ske04p'}. In fact, $E$ may be viewed as \it{Morita equivalence} from $\sK(E)$ to $\cB$. The inverse operation, conjugation of $E_t$ with $E$, gives $E\odot E_t\odot E^*=\cls\sK(E)F_t\sK(E)$. The left action of $\sK(E)$ is nondegenerate by strictness. But, there is no reason why the inner product of $F_t$ should assume values only in $\sK(E)$. In fact, if $T$ was a Markov semigroup, then each $F_t$ has a unit vector and $\sK(E)$ is rarely unital.

When $T$ is an \hl{\nbd{E}semigroup} (that is, the maps $T_t$ are endomorphisms of $\sB^a(E)$), or even an \hl{\nbd{E_0}semigroup} (that is, the endomorphisms are also unital), then the product system is the one associated in Skeide \cite{Ske04p'} with a strict \nbd{E}semigroup. (See Skeide \cite{Ske02} for the first construction for \nbd{E_0}semigroup when $E$ has a unit vector, and Bhat and Lindsay \cite{BhLi05} for an \nbd{E}semigroup under the same hypothesis.) When $E=H$ is a Hilbert space, we recover Bhat's construction \cite{Bha96} of the Arveson system of a normal \nbd{E_0}semigroup on $\sB(H)$. When $T$ is a normal CP-semigroup on $\sB(G)$, we obtain a direct construction of its Arveson system. (In Bhat \cite{Bha96}, it is constructed via a so-called \it{minimal dilation} of $T$ to an \nbd{E}semigroup on $\sB(H)$ as the Arveson system of that \nbd{E}semigroup.)
\erem

The product system $E^\odot$ has no relation with $T$ as direct as the GNS-system $F^\odot$. (There is no unit for $E^\odot$ that would allow to uncover the CP-semigroup $T$. In fact, if $E=H$ is a Hilbert space, then it is known that $E^\odot$ can be unitless. In \cite[Theorem 3.4]{BLS02p} we have shown that the product system of the minimal dilation is $E^\odot$.) But the following theorem shows that spatiality is preserved. For von Neumann modules this is a new result in the classification of product systems up to Morita equivalence; see Remark \ref{vNMerem}.

\bthm\label{BaEspatthm}
Let $E$ be a full Hilbert module over a unital \nbd{C^*}algebra $\cB$. Suppose $F^\odot$ is a product system of correspondences $F_t$ over $\sB^a(E)$ with strict left actions. Define the product system $E^\odot$ as above by setting $E_t=E^*\odot F_t\odot E$.

If $F^\odot$ is spatial, then so is $E^\odot$. More precisely, if $\Om^\odot$ is the central unital reference unit of $F^\odot$, then by $i_t\colon\AB{x,y}\mapsto x^*\odot\Om_t\odot y\in E_t$ for each $t\in\R_+$, we define an injective morphism from the trivial product system $\bfam{\cB}_{t\in\R_+}$ into $E^\odot$. In particular, the image of the central unital unit $\bfam{\U}_{t\in\R_+}$ is a central unital unit $\om^\odot$ for $E^\odot$.
\ethm

\proof
For each $t\in\R_+$, the map $i_t$ is an isometry. Indeed,
\beqn{
\AB{x^*\odot\Om_t\odot y,x'^*\odot\Om_t\odot y'}
~=~
\AB{\Om_t\odot y,xx'^*\Om_t\odot y'}
~=~
\AB{\Om_t\odot y,\Om_t\odot xx'^*y'}
~=~
\AB{y,x}\AB{x',y'}.
}\eeqn
Clearly, $i_t$ is bilinear. Since $\cB$ is unital and $E$ is full, by \cite[Lemma 3.2]{Ske04p'} there exist $n\in\N$ and $x_1,\ldots,x_n\in E$ such that $\sum_{i=1}^n\AB{x_i,x_i}=\U$. So, for $\om_t:=i_t(\U)=\sum_{i=1}^nx_i^*\odot\Om_t\odot x_i$ we find
\bmun{
\om_s\odot\om_t
~=~
\sum_{i,j=1}^nx_i^*\odot\Om_s\odot x_i\odot x_j^*\odot\Om_t\odot x_j
~=~
\sum_{i,j=1}^nx_i^*\odot\Om_s\odot x_ix_j^*\Om_t\odot x_j
\\
~=~
\sum_{i,j=1}^nx_i^*\odot\Om_s\odot\Om_t\odot x_ix_j^*x_j
~=~
\sum_{i=1}^nx_i^*\odot\Om_{s+t}\odot x_i
~=~
\om_{s+t},
}\emun
so that the $\om_t$ form a unit $\om^\odot$ for $E^\odot$. By bilinearity of $i_t$, the unit $\om^\odot$ is unital and central. By the same reason, the $i_t$ form a morphism, that is, $i_s\odot i_t=i_{s+t}$.\qed

\brem
A similar result is true for von Neumann modules. Just that one has to refer to \cite[Lemma 4.2]{Ske04p'} and the sum is no longer finite and, in general, only strongly convergent.
\erem

\brem\label{vNMerem}
Note that the converse statement may fail. The simplest reason is that the correspondences $F_t$ of $F^\odot$ need not allow for unit vectors. More concretely, $F_t$ need not be strictly complete. (Otherwise, choose a bounded approximate unit for $\sK(E)$ that consists of finite-rank operators $\sum_{i=1}^nx_iy_i^*$. Then the corresponding net of elements $\sum_{i=1}^nx_i\odot\om_t\odot y_i^*$ is strictly Cauchy in $F_t$. If the limits $\Om_t$ exist, then they form a unital central unit for $F^\odot$.) For von Neumann modules also the converse statement is true: $E^\odot$ is spatial if and only if $F^\odot$ is spatial. (Von Neumann modules are not only strictly complete but even \nbd{\sigma}weakly.) In the sense of Morita equivalence of product systems of von Neumann correspondences \cite{Ske04p'}, one may rephrase as follows: Morita equivalence of product systems of von Neumann correspondences preserves spatiality.
\erem

Now, since we know what is the product system $E^\odot$ of correspondences over $\cB$ of a strict CP-semigroup $T$ on $\sB^a(E)$ for some full Hilbert \nbd{\cB}module $E$, and since we know that spatiality of $T$ is reflected by spatiality of (some spatial extension) $E^\odot$, we can ask whether there possibly is a Powers sum for spatial CP-semigroups such that the sum operation is reflected by the product operation of their spatial product systems of correspondences over $\cB$. For \nbd{E_0}semigroups we proved the affirmative answer in \cite{BLS02p}. For spatial CP-semigroups the result is new. The proof also simplifies the proof of \cite{BLS02p}.

We start with a simple consequence of Observation \ref{gensetob}.

\blem
Let $E$ be a full Hilbert module over a unital \nbd{C^*}algebra $\cB$. Let $T$ be a spatial strict CP-semigroup on $\sB^a(E)$ and choose a unit $S$ implemented by a semigroup $c$ in $\sB^a(E)$. Denote by $F^\odot$ the spatial extension of the GNS-system associated with that unit as in the proof of Theorem \ref{spatCPDthm} (considering $T$ a CPD-semigroup on a one-point set), so that $F^\odot$ is generated by the unit $\zeta^\odot$ that gives back $T_t=\AB{\zeta_t,\bullet\zeta_t}$ and by the central unital reference unit $\Om^\odot$. Denote by $E^\odot$ and $\om^\odot$  product system and central unit as in Theorem \ref{BaEspatthm}.

Then the product system $E^\odot$ is generated by elements of the form $\om_t$ and $x^*\odot\zeta_t\odot y$ in the sense that
\beqn{
E_t
~=~
\cls\BCB{z_n\odot\ldots\odot z_1\colon n\in\N,t_1+\ldots+t_n=t,z_i=\om_{t_i}\text{~or~}z_i=x^*\odot\zeta_{t_i}\odot y}.
}\eeqn
\elem

We omit the obvious proof. Note, howewer, that
\baln{
\AB{x^*\odot\zeta_t\odot y,x'^*\odot\zeta_t\odot y'}
&
~=~
\AB{y,T_t(xx'^*)y'},
&
\AB{x^*\odot\zeta_t\odot y,x'^*\odot\Om_t\odot y'}
&
~=~
\AB{y,c_t^*x}\AB{x',y},
}\ealn
so that $\AB{x^*\odot\zeta_t\odot y,\om_t}=\AB{y,c_t^*x}$. Note, too, that the pair $(E^\odot,\om^\odot)$ is determined by these properties up to spatial isomorphism.

\bthm\label{PsCPthm}
For $i=1,2$, let $T^i$ be spatial strict CP-semigroup on $\sB^a(E^i)$ for full Hilbert modules $E^i$ over a unital \nbd{C^*}algebra $\cB$. Choose units $S^i$ for $T^i$ implemented by semigroups $c^i$ in $\sB^a(E^i)$. Then, by setting
\beqn{
(T^1\boxplus T^2)_t\fMatrix{a_{11}&a_{12}\\a_{21}&a_{22}}
~:=~
\fMatrix{T^1_t(a_{11})&{c^1_t}^*a_{12}c^2_t\\{c^2_t}^*a_{21}c^1_t&T^2_t(a_{22})},
}\eeqn
we define a spatial CP-semigroup on $\sB^a(E^1\oplus E^2)$ with a unit $S^1\boxplus S^2$ implemented by $c=c^1\oplus c^2$, the \hl{Powers sum} $(T^1\boxplus T^2,S^1\boxplus S^2)$ of $(T^1,S^1)$ and $(T^2,S^2)$. Moreover, if $(E^\odot,\om^\odot)$ denotes the spatial product system of correspondences over $\cB$ associated with $(T^1\boxplus T^2,S^1\boxplus S^2)$, and if $({E^i}^\odot,{\om^i}^\odot)$ denote those associated with $(T^i,S^i)$, then $(E^\odot,\om^\odot)=(E^1\circledcirc E^2,\om^\odot)$.
\ethm

\proof
The proof that $T^1\boxplus T^2$ is a CP-semigroup and that it is spatial with unit $c=c^1\oplus c^2$ in $\sB^a(E^1\oplus E^2)$, is similar to that of Theorem \ref{PsCPDthm}.

To proof that $E^\odot$ is product of ${E^1}^\odot$ and ${E^2}^\odot$, we observe that by the lemma ${E^i}^\odot$ is generated by expressions $\om^i_t$ and ${x^i}^*\odot\zeta^i\odot y^i$. Therefore the product $E^1\circledcirc E^2$ is generated by expressions $\om_t$, ${x^1}^*\odot\zeta^1\odot y^1$, and ${x^2}^*\odot\zeta^2\odot y^2$, where the only yet unspecified inner product is
\beqn{
\AB{{x^1}^*\odot\zeta^1\odot y^1,{x^2}^*\odot\zeta^2\odot y^2}
~=~
\AB{{x^1}^*\odot\zeta^1\odot y^1,\om^1_t}\AB{\om^2_t,{x^2}^*\odot\zeta^2\odot y^2}
~=~
\AB{y^1,{c^1_t}^*x^1}\AB{{c^2_t}^*x^2,y^2}.
}\eeqn
On the other hand, $E^\odot$ is generated by expressions $\om_t$ and $\raisebox{1pt}{\tMatrix{x^1\\x^2}}^*\odot\zeta_t\odot\raisebox{1pt}{\tMatrix{y^1\\y^2}}$. By calculating the norm, one easily verifies that $\raisebox{1pt}{\tMatrix{x^1\\0}}^*\odot\zeta_t\odot\raisebox{1pt}{\tMatrix{0\\y^2}}=\raisebox{1pt}{\tMatrix{0\\x^2}}^*\odot\zeta_t\odot\raisebox{1pt}{\tMatrix{y^1\\0}}=0$, while
\baln{
\BAB{\raisebox{1pt}{\sMatrix{x^1\\0}}^*\odot\zeta_t\odot\raisebox{1pt}{\sMatrix{y^1\\0}}\,\,,\,\raisebox{1pt}{\sMatrix{x'^1\\0}}^*\odot\zeta_t\odot\raisebox{1pt}{\sMatrix{y'^1\\0}}}
&
~=~
\AB{{x^1}^*\odot\zeta^1_t\odot y^1,{x'^1}^*\odot\zeta^1_t\odot y'^1},
\\
\BAB{\raisebox{1pt}{\sMatrix{0\\x^2}}^*\odot\zeta_t\odot\raisebox{1pt}{\sMatrix{0\\y^2}}\,\,,\,\raisebox{1pt}{\sMatrix{0\\x'^2}}^*\odot\zeta_t\odot\raisebox{1pt}{\sMatrix{0\\y'^2}}}
&
~=~
\AB{{x^2}^*\odot\zeta^2_t\odot y^2,{x'^2}^*\odot\zeta^2_t\odot y'^2},
\\
\BAB{\raisebox{1pt}{\sMatrix{x^1\\0}}^*\odot\zeta_t\odot\raisebox{1pt}{\sMatrix{y^1\\0}}\,\,,\,\raisebox{1pt}{\sMatrix{0\\x^2}}^*\odot\zeta_t\odot\raisebox{1pt}{\sMatrix{0\\y^2}}}
&
~=~
\AB{{x^1}^*\odot\zeta^1_t\odot y^1,\om^1_t}\AB{\om^2_t,{x^2}^*\odot\zeta^2_t\odot y^2}.
}\ealn
So, $E^\odot$ is isomorphic to $E^1\circledcirc E^2$, via
\baln{
\om_t
&
~\longmapsto~
\om_t,
&
\raisebox{1pt}{\sMatrix{x^1\\0}}^*\odot\zeta_t\odot\raisebox{1pt}{\sMatrix{y^1\\0}}
&
~\longmapsto~
{x^1}^*\odot\zeta^1_t\odot y^1,
&
\raisebox{1pt}{\sMatrix{0\\x^2}}^*\odot\zeta_t\odot\raisebox{1pt}{\sMatrix{0\\y^2}}
&
~\longmapsto~
{x^2}^*\odot\zeta^2_t\odot y^2.\qedsymbol}\ealn
\noqed

\vspace{-2ex}
\brem
The algebras $\sB^a(E^1)$ and $\sB^a(E^2)$ have the property that they may be considered as subalgebras of the \it{matrix algebra} (see Skeide \cite{Ske00b} for details about matrix algebras) $\sB^a(E^1\oplus E^2)=\raisebox{2pt}{\tMatrix{\sB^a(E^1)&\sB^a(E^2,E^1)\\\sB^a(E^1,E^2)&\sB^a(E^2)}}$. The interesting property is that both $\sB^a(E^1)$ and $\sB^a(E^2)$ are generated from products of the off-diagonal entries in the strict topology. (Note that this may fail, if $E^1$ and $E^2$ are not both full, up to the point where $\sB^a(E^1,E^2)=\zero$.) In the case of von Neumann modules that means that $\sB^a(E^1)$ and $\sB^a(E^2)$ are Morita equivalent as von Neumann algebras. (In fact, both are Morita equivalent as von Neumann algebra to $\cB$.) For \nbd{C^*}modules one might say, they are \it{strictly} Morita equivalent. (We do not know whether there exists a systematic study of Morita equivalence for multiplier algebras. This would be adequate to our purposes, as the multiplier algebra of $\sK(E)$ is $\sB^a(E)$. We met already several times, in \cite{Ske04p',BLS02p}, situations where we had to develop at least parts of such a theory.)
\erem

\brem
The case when $T^i$ are \nbd{E_0}semigroups has been discussed in \cite{BLS02p}. The proof here, restricted to that case, differs considerably from that in \cite{BLS02p} and, actually, appears simpler. The case of \nbd{E_0}semigroups acting on $\sB(H)$s, is the one proposed by Powers 2002 in Mount Holyoke; see also Example \ref{SchurCPex}. The Arveson system of the Powers sum in that case has been identified as product in \cite{Ske03c}. The case of CP-semigroups acting on $\sB(H)$s has been discussed in \cite{Pow04} with a much less general notion of spatiality for CP-semigroups. \cite{Pow04} also does not identify the Arveson system of the sum as product of spatial Arveson systems. But, he proves that it need not be isomorphic to the tensor product (available only for Arveson systems).
\erem

\brem
Of course, like in Remark \ref{Afamrem}, also here all statements remain true for families of spatial CP-semigroups and the spatial extensions of their GNS-systems.
\erem

\section{CPD-semigroups \it{versus} Schur CP-semigroups}\label{vsSEC}

In Example \ref{SchurCPex} we pointed out that \nbd{\sB(\cB)}valued CPD-semigroups on a finite set $S$ (with cardinality $n$, say) are in one-to-one correspondence with Schur CP-semigroups on $M_n(\cB)$ and that this one-to-one correspondence behaves well with respect to the respective Powers sums. After Theorem \ref{PsCPthm}, we can say that this one-to-one correspondence also behaves well with respect to the products of the respective spatial extensions of the product systems of correspondences over $\cB$. (They simply coincide.)

In this section we wish to see in how far we can generalize that one-to-one correspondence to arbitrary sets $S$. The idea in Example \ref{SchurCPex} was to let act the semigroups $\eT^{s,s'}$ on the matrix elements $a_{s,s'}$ of a finite \nbd{\#S\times\#S}matrix with entries in $\cB$. We simply try now to do the same with \nbd{\#S\times\#S}matrices for a set $S$ of arbitrary cardinality $\#S$.

Of course, the matrices should continue to form a \nbd{C^*}algebra, so we cannot allow arbitrary matrices. A canonical candidate is the \nbd{C^*}algebra $\sB^a(\cB^S)$ where $\cB^S$ is the \nbd{\#S}dimensional \hl{column space} space of $\cB$. $\cB^S$ consists of all families $B=\bfam{b_s}_{s\in S}$ such that the net $\sum_{s\in S'}b_s^*b_s$ converges over the finite subsets $S'$ of $S$. The inner product is $\AB{B,B'}:=\sum_{s\in S}b_s^*b'_s$.

Let $e_s:=\bfam{\delta_{ss'}\U}_{s'\in S}$. The elements $e_s$ form an \it{orthonormal basis} of $\cB^S$ in the obvious way: $\sum_{s\in S}e_ie_i^*=\id_{\cB^S}$ strongly and, therefore, \nbd{*}strongly in $\sB^a(\cB^S)$ over the finite subsets of $S$, and since the approximating net is bounded by $\U$, also strictly. It follows that an arbitrary element $A\in\sB^a(\cB^S)$ can be written as
\beqn{
A
~=~
\sum_{s,s'\in S}e_sa_{s,s'}e_{s'}^*,
}\eeqn
where $a_{s,s'}:=\AB{e_s,Ae_{s'}}\in\cB$. (We resist the temptation to denote $\sB^a(\cB^S)$ as $M_S(\cB)$, because the latter, usually, rather refers to $\sK(\cB^S)$.)

A \hl{Schur CP-map} on $\sB^a(\cB^S)$ is a CP-map $T$ on $\sB^a(\cB^S)$ such that
\beqn{
T(e_se_s^*Ae_{s'}e_{s'}^*)
~=~
e_se_s^*T(A)e_{s'}e_{s'}^*
}\eeqn
for all $A\in\sB^a(\cB^S)$ and all $s,s'\in S$. Without the simple proof we state the following:

\bprop
A Schur CP-map necessarily leaves $\sK(\cB^S)$ invariant and is strict.
\eprop

Obviously, if $T$ is a Schur CP-map, then $\eK^{s,s'}:=\AB{e_s,T(e_s\bullet e_{s'}^*)e_{s'}}$ defines a \nbd{\sB(\cB)}valued CPD-kernel $\eK$ on $S$. Moreover, $T$ can be recovered from $\eK$ as
\beq{\label{CPCPD}
T(e_sbe_{s'}^*)
~=~
e_s\eK^{s,s'}(b)e_{s'}^*.
}\eeq
However, not all CPD-kernels give rise to Schur CP-map in that way. A CPD-kernel is \hl{bounded}, if there is a constant $M$ such that $\snorm{\eK^{s,s'}}\le M$ for all $s,s'\in S$.

\bprop
Let $\eK$ be a \nbd{\sB(\cB)}valued CPD-kernel $\eK$ on $S$. Then $\eK$ gives rise to a (unique) Schur CP-map on $\sB^a(\cB^S)$ fulfilling \eqref{CPCPD} if and only if $\eK$ is bounded.
\eprop

\proof
Clearly, a kernel $\eK$ fulfilling \eqref{CPCPD} for some CP-map $T$, is bounded by $M=\norm{T}$. So, for the other direction let us suppose that $\eK$ is bounded (by $M$, say). Instead of showing directly that under this condition the map defined by \eqref{CPCPD} on finite matrices (that it, operators $A\in\sB^a(\cB^S)$ with only finitely many matrix entries $a_{s,s'}$ different from $0$) extends suitably to a CP-map $T$, we shall construct a candidate for the GNS-construction of $T$.

Let $(E,\bfam{\xi^s}_{s\in S})$ denote the Kolmogorov decomposition for $\eK$. Define $F:=\cB^S\odot E\odot\cB_S$, where $\cB_S:={\cB^S}^*$, the \nbd{\#S}dimensional \hl{row space} of $\cB$. Recall that an element $y\in F$ may be interpreted as a map $B\mapsto y\odot B$ from $\cB^S$ to $F\odot\cB^S=\cB^S\odot E$. We claim that the sum $\sum_{s\in S}e_s\odot\xi^s\odot e_s^*$ converges \nbd{*}strongly in $\sB^a(\cB^S,\cB^S\odot E)$ to an operator $Z$. Once we have convergence, it is clear that the CP-map $T(A):=Z^*(A\odot\id_E)Z$ fulfills \eqref{CPCPD}.

Let $B=\bfam{b_s}_{s\in S}\in\cB^S$. Then for every finite subset $S'\subset S$ we have
\beqn{
\Babs{\sum_{s\in S'}(e_s\odot\xi^s\odot e_s^*)B}^2
~=~
\Babs{\sum_{s\in S'}e_s\odot\xi^sb_s}^2
~=~
\sum_{s\in S'}b_s^*\AB{\xi^s,\xi^s}b_s
~\le~
M\sum_{s\in S'}b_s^*b_s.
}\eeqn
From this two things follow. Firstly, $\sum_{s\in S'}(e_s\odot\xi^s\odot e_s^*)B$ is a Cauchy net in $\cB^S\odot E$. Secondly, the net $\sum_{s\in S'}e_s\odot\xi^s\odot e_s^*$ is bounded by $\sqrt{M}$. From boundedness it follows that strong convergence of the adjoint net may be checked on the total subset $e_s\odot x$ of $\cB^S\odot E$. But this is clear because $\sum_{s\in S'}(e_s\odot\xi^s\odot e_s^*)^*(e_{s'}\odot x)=e_{s'}\AB{\xi^{s'},x}$ if $s'\in S'$ and $0$ otherwise.\qed

\bcor
Fix a set $S$ and a unital \nbd{C^*}algebra $\cB$. Then the formula \eqref{CPCPD}, when applied to all members of a semigroup, establishes a one-to-one correspondence between \hl{pointwise bounded} \nbd{\sB(\cB)}valued CPD-semigroups $\eT$ on $S$ (that is, each $\eT_t$ is bounded) and Schur CP-semigroups $T$ on $\cB^S$. Moreover, $\eT$ and $T$ have the same product systems of correspondences over $\cB$.
\ecor

\proof
We do not prove the only still open statement about the product systems. (It follows from the observation that the unit $\zeta^\odot$ of the GNS-system of $T$ is a \it{diagonal matrix} with the unit ${\xi^s}^\odot$ of the GNS-system of $\eT$ as \nbd{s,s}entry; see also \cite[Appendix B]{BBLS04}.)\qed

\bthm
$\eT$ is spatial if and only if $T$ is spatial. In that case, the one-to-one correspondence respects also units and the spatial extensions of the GNS-systems based on them. Therefore, it must respect also Powers sum (and obviously products of the product systems, because the product systems coincide, anyway).
\ethm

We omit the obvious proof also here.

\brem
The formulation is for \nbd{C^*}algebras and modules. A similar correspondence has been proved for von Neumann algebras and modules in \cite[Appendix B]{BBLS04}, however, only for uniformly continuous semigroups, and without paying attention to one-to-one aspect and the related notion of Schur CP-semigroup. (Recall that uniform continuity is automatic only for spatial semigroups. The statements that are valid also in the nonspatial case, do not require any continuity in time.) Anyway, all statements remain true also in the von Neumann case, some of them simpler, because no spatial extension is needed.
\erem

\brem
Note that for pointwise bounded CPD-semigroups, the results in Section \ref{CPDPsSEC} maybe obtained from those in Section \ref{CPPsSEC} via the one-to-one correspondence. (The only exception is the spatial extension of the GNS-system of a spatial CP-semigroup. But this can easily be added to Section \ref{CPPsSEC}, to make it independent of Section \ref{CPDPsSEC}.) The not pointwise bounded case can also be reduced to the pointwise bounded case, rescaling the CPD-semigroup with scalar semigroups. But this discussion is somewhat cumbersome and not at all instructive. We prefer to leave Section \ref{CPDPsSEC} as a separate one, which in its general form is not included in Section \ref{CPPsSEC}.
\erem



\setlength{\baselineskip}{2.5ex}


\begin{thebibliography}{BBLS04}

\bibitem[Arv89]{Arv89a}
W.~Arveson, \emph{{Continuous analogues of Fock space III: Singular states}},
  J.\ Operator Theory \textbf{22} (1989), 165--205.

\bibitem[Arv97]{Arv97a}
\bysame, \emph{{The index of a quantum dynamical semigroup}}, J.\ Funct.\ Anal.
  \textbf{146} (1997), 557--588.

\bibitem[BBLS04]{BBLS04}
S.D. Barreto, B.V.R. Bhat, V.~Liebscher, and M.~Skeide, \emph{{Type I product
  systems of Hilbert modules}}, J.\ Funct.\ Anal. \textbf{212} (2004),
  121--181, (Preprint, Cottbus 2001).

\bibitem[Bha96]{Bha96}
B.V.R. Bhat, \emph{{An index theory for quantum dynamical semigroups}}, Trans.\
  Amer.\ Math.\ Soc. \textbf{348} (1996), 561--583.

\bibitem[BL05]{BhLi05}
B.V.R. Bhat and J.M. Lindsay, \emph{{Regular quantum stochastic cocycles have
  exponential product systems}}, Quantum Probability and Infinite Dimensional
  Analysis --- From Foundations to Applications (M.~Sch\"urmann and U.~Franz,
  eds.), Quantum Probability and White Noise Analysis, no. XVIII, World
  Scientific, 2005, pp.~126--140.

\bibitem[BLS07]{BLS02p}
B.V.R. Bhat, V.~Liebscher, and M.~Skeide, \emph{{A problem of Powers and the
  product of spatial product systems}}, ar\-Xiv: 0801.0042v1, 2007.

\bibitem[BLS08]{BLS08p}
\bysame, \emph{{Subsystems of Fock need not be Fock: Spatial CP-semigroups}},
  ar\-Xiv: 0804.2169v1, 2008.

\bibitem[BS00]{BhSk00}
B.V.R. Bhat and M.~Skeide, \emph{{Tensor product systems of Hilbert modules and
  dilations of completely positive semigroups}}, Infin.\ Dimens.\ Anal.\
  Quantum Probab.\ Relat.\ Top. \textbf{3} (2000), 519--575, (Rome,
  Volterra-Pre\-print 1999/0370).

\bibitem[Lie03]{Lie00p1}
V.~Liebscher, \emph{{Random sets and invariants for (type II) continuous tensor
  product systems of Hilbert spaces}}, Pre\-print, ar\-Xiv: math.PR/0306365,
  2003, To appear in Mem.\ Amer.\ Math.\ Soc.

\bibitem[Pas73]{Pas73}
W.L. Paschke, \emph{{Inner product modules over $B^*$--algebras}}, Trans.\
  Amer.\ Math.\ Soc. \textbf{182} (1973), 443--468.

\bibitem[Pow04]{Pow04}
R.T. Powers, \emph{{Addition of spatial $E_0$--semigroups}}, Operator algebras,
  quantization, and noncommutative geometry, Contemporary Mathematics, no. 365,
  American Mathematical Society, 2004, pp.~281--298.

\bibitem[PS72]{PaSchm72}
K.R. Parthasarathy and K.~Schmidt, \emph{{Positive definite kernels, continuous
  tensor products, and central limit theorems of probability theory}}, Lect.\
  Notes Math., no. 272, Springer, 1972.

\bibitem[Ske00]{Ske00b}
M.~Skeide, \emph{{Generalized matrix $C^*$--algebras and representations of
  Hilbert modules}}, Mathematical Proceedings of the Royal Irish Academy
  \textbf{100A} (2000), 11--38, (Cott\-bus, Rei\-he Mathe\-ma\-tik 1997/M-13).

\bibitem[Ske02]{Ske02}
\bysame, \emph{{Dilations, product systems and weak dilations}}, Math.\ Notes
  \textbf{71} (2002), 914--923.

\bibitem[Ske03a]{Ske03c}
\bysame, \emph{{Commutants of von Neumann modules, representations of
  $\sB^a(E)$ and other topics related to product systems of Hilbert modules}},
  Advances in quantum dynamics (G.L. Price, B~.M. Baker, P.E.T. Jorgensen, and
  P.S. Muhly, eds.), Contemporary Mathematics, no. 335, American Mathematical
  Society, 2003, (Preprint, Cottbus 2002, ar\-Xiv: math.OA/0308231),
  pp.~253--262.

\bibitem[Ske03b]{Ske03b}
\bysame, \emph{{Dilation theory and continuous tensor product systems of
  Hilbert modules}}, Quantum Probability and Infinite Dimensional Analysis
  (W.~Freudenberg, ed.), Quantum Probability and White Noise Analysis, no.~XV,
  World Scientific, 2003, Preprint, Cottbus 2001, pp.~215--242.

\bibitem[Ske04]{Ske04p'}
\bysame, \emph{{Unit vectors, Morita equivalence and endomorphisms}},
  Pre\-print, ar\-Xiv: math.OA/0412231v5 (Version 5), 2004, To appear in Publ.\
  Res.\ Inst.\ Math.\ Sci.

\bibitem[Ske06]{Ske06d}
\bysame, \emph{{The index of (white) noises and their product systems}},
  Infin.\ Dimens.\ Anal.\ Quantum Probab.\ Relat.\ Top. \textbf{9} (2006),
  617--655, (Rome, Volterra-Pre\-print 2001/0458, ar\-Xiv: math.OA/0601228).

\bibitem[Ske08a]{Ske08p1}
\bysame, \emph{{Classification of $E_0$--semigroups by product systems}},
  Pre\-print, in preparation, 2008.

\bibitem[Ske08b]{Ske08}
\bysame, \emph{{Isometric dilations of representations of product systems via
  commutants}}, Int.\ J.\ Math. \textbf{19} (2008), 521--539, (ar\-Xiv:
  math.OA/0602459).

\end{thebibliography}
\newcommand{\Swap}[2]{#2#1}\newcommand{\Sort}[1]{}
\providecommand{\bysame}{\leavevmode\hbox to3em{\hrulefill}\thinspace}
\providecommand{\MR}{\relax\ifhmode\unskip\space\fi MR }
\providecommand{\MRhref}[2]{%
  \href{http://www.ams.org/mathscinet-getitem?mr=#1}{#2}
}
\providecommand{\href}[2]{#2}


\end{document}